\documentclass[12pt]{article}

\usepackage{graphics}
\usepackage[hang,flushmargin]{footmisc}
\usepackage{sans}
\usepackage{amsmath}
\usepackage{amssymb} 
\usepackage{fullpage}
\def\halfneg{{\setbox1=\hbox{$\!$}\hspace{0.5\wd1}}}

\def\deRham{de~Rham}

\def\dee#1{{\mathit{d}{#1}}}

\def\integralOf#1#2#3#4#5{{\displaystyle{#5}_{#3}^{#4}{#1}\,{#2}}}
\def\singleIntegral#1#2#3#4{{\integralOf{#1}{#2}{\!\!\!{#3}}{#4}{\int}}}

\def\Hom{\operatorname{Hom}}

\def\Totto{\longrightarrow}

\begin{document}

\hrule \textbf{\Large Georges \deRham\ 1903--1990}\\ \hrule

\renewcommand{\thefootnote}{\fnsymbol{footnote}}
{\scriptsize 
 This is an English translation of the obituary notice by Beno Eckmann,
 appearing under the same title in \textit{\textsf{Elemente der 
 Mathematik}\/} \textbf{47}(3) 118--122 (1992) (in 
 German).\footnote{\textit{\textsf{Translator's note:}\/} Prepared by
 Neil Ching (2008) and distributed with permission of the publisher.
 I am indebted to Professor Urs Stammbach of ETH Z\"urich for his
 assistance.} Original \copyright\ 1992 Birkh\"auser Verlag and
 available at \texttt{<http://eudml.org/doc/141534>}.}
\renewcommand{\thefootnote}{\arabic{footnote}}

\vspace{\parskip}
\begin{center}
\scalebox{0.64}{\includegraphics{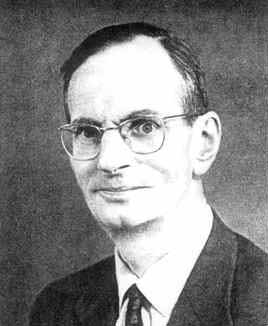}}
\end{center}
\vspace{\parskip}

The news of the passing of Georges \deRham\ on 8 October 1990 has by
now surely reached and moved all in the mathematical community.  This
journal would like to bid him farewell with the following
communication. He was one of the important figures of mathematics in
our century, his name and work belong to its enduring legacy, and the
influence of his ideas on its development has by no means been
exhausted yet.

He was close to \textit{\textsf{Elemente der Mathematik}\/}, which is
broadly aimed at both instructors and students, and indeed for many
years was among its contributors. Everything that was connected with
teaching lay as close to his heart as research did. He was in the
habit of saying, ``Teaching, the conveying of essentials, making the
beautiful intelligible and evident, that is what gives me joy; and
instruction is always accompanied by interpretation.''

He understood how to enlighten, in an unassuming yet memorable
fashion, students at all levels about mathematics; perhaps they could
unconsciously sense what a great mathematical powerhouse was at work
here. The closer to him one became, whether as a novice or as a
colleague, the more one was impressed by his personality: by
his---there is no other way to say it---refined, yet not distant,
bearing; by his charm, which came from the heart; by his unyielding
straightforwardness and intensity; by his warmhearted friendship,
loyalty, and readiness to help.

\vspace{\parskip}
\textbf{\large Biographical Sketch}

Georges \deRham\ came from a prominent family of the Swiss canton of
Vaud.  He grew up in Lausanne and studied there at the university,
receiving his \textit{\textsf{licence}\/}. In Paris, where he spent
several years, he came into contact with the work of Henri Poincar\'e
and \'Elie Cartan. Here he was led to focus his research on the
mathematical area to which his 1931 doctoral thesis is devoted, a
piece of work that immediately made his name world-famous. He
qualified as university lecturer at Lausanne in 1932, and from 1936
until he took emeritus status he was professor at the universities of
Lausanne and Geneva; he also maintained his contacts with Paris
throughout. Although many centers of mathematics all over the world
strove to lure him, he never left Romandy, and he succeeded in
creating a world-class school of mathematics there. He encouraged
young talent in a selfless, collegial, and fatherly way.  With great
vigor, he invested much time and effort providing services to the
larger scientific world, as president of the Swiss Mathematical
Society during the difficult time in 1944/45, as editor of
\textit{\textsf{Commentarii Mathematici Helvetici}\/}, and as a member
of the research council of the Swiss National Science Foundation. From
1963 till 1966 the International Mathematical Union entrusted him with
its highest office, that of president, and here, thanks to his
reliable and universally trusted leadership, he succeeded in promoting
worldwide scientific cooperation.

His friends knew that he was accomplished not only as a mathematician,
but as an alpinist as well. In this endeavor too he felt himself drawn
toward difficult challenges. Thanks to first ascents in the Alps, an
excellently written book, and lectures, \deRham's name was a household
word in mountainclimbing circles. After the war, when the borders were
once again open, he recounted with pride that his first invited
lecture abroad was not for mathematicians, but for the English
``Rucksack Club.''

It is hardly astonishing that many honors, local, national, and
international, were bestowed upon \deRham. This is not the place to list
them all, and we mention only three here: honorary member of the Swiss
Mathematical Society, 1960; honorary doctorate, ETH Zurich, 1961; and the
federally awarded ``Marcel Benoist Prize,'' 1966.

\vspace{\parskip}
\textbf{\large De~Rham's Theorems}

De~Rham's mathematical work lies in the areas of global analysis,
differential geometry, and topology. It is striking for its
astonishing synthesis of the most disparate of methods: analytical,
geometric, combinatorial-topological, and algebraic. It was typical of
him that he worked only within the scope of a few specific topics, but
impressed upon them the stamp of his powers of penetration. They are
as follows: differential forms on manifolds; combinatorial invariants
of cell complexes; the definitive version of the Hodge theory of
harmonic differential forms; the decomposition of Riemannian manifolds
into products of indecomposable submanifolds; and ``currents,'' a
synthesis, making use of Schwartz distributions, of differential forms
and the chains of algebraic topology.

Here we restrict ourselves to summarizing his contribution to the
first subject---it is the most well-known and the one that is most
often applied.

The famous \textit{\textsf{\deRham's theorem}\/}, announced in 1928
and worked out in his \textit{\textsf{th\`ese}\/} (Paris, 1931), is
today usually formulated, in the terminology of cohomology theory, as
an isomorphism:
\begin{equation}
 R^p(M)\ \cong\ H^p(M)
   \textrm{.}
\end{equation}
Let us briefly recall the meaning of these symbols.

We consider differential forms of degree $p$, $0\leq p\leq n$, on a closed
differentiable manifold $M$ of dimension $n$; they make up a real vector
space $\Omega^p$. The exterior differential $\dee{\alpha^p}$ of $\alpha^p
\in\Omega^p$ provides a linear mapping
\[ \dee{}:\Omega^p\Totto\Omega^{p+1}\textrm{,}
\]
and we have $\dee{}\dee{}=0$; we set $\Omega^{-1}=\Omega^{n+1}=0$. With
this, the image $\dee{\Omega^{p-1}}$ lies in the kernel $\mathring\Omega^p$
of $\dee{}$, and the ``cohomology groups''
\[ \mathring\Omega^p/\dee{\Omega^{p-1}}=R^p(M)
\]
are the \textit{\textsf{\deRham\ groups}\/} of $M$ (which are real
vector spaces). A form $\alpha^p\in\mathring\Omega^p$, which thus has
the differential $\dee{\alpha^p}=0$, is said to be
\textit{\textsf{closed}\/}.

Further suppose that $M$ is subdivided into cells, and let $C_p$ be
the vector space of $p$-chains, that is, the vector space of linear
combinations, with real coefficients, of $p$-dimensional cells. The
boundary $\partial\sigma_p$ of a $p$-cell $\sigma_p$ is a
$(p-1)$-chain (with coefficients $\pm1$), and through linearity we
obtain a mapping
\[ \partial:C_p\Totto C_{p-1}
\]
such that $\partial\partial=0$; we set $C_{-1}=C_{n+1}=0$. The vector
space dual to $C_p$, that is, the space of the real linear forms
defined on $p$-chains, will be denoted by $C^p$. For these, we have
the ``coboundary'' operator
\[ \delta:C^p\Totto C^{p+1}\textrm{,}
\]
which is dual to $\partial$ and defined for $f^p\in C^p$ by $\delta
f^p(c_{p+1})=f^p(\partial c_{p+1})$, where $c_{p+1}\in C_{p+1}$; we
have $\delta\delta=0$. Hence $\delta C^{p-1}$ lies in the kernel
$\mathring{C}^p$ of $\delta$, and one obtains the cohomology groups
\[ H^p(M)\ =\ \mathring{C}^p/\delta C^{p-1}
\]
(which again are real vector spaces). These are topological invariants
of $M$, and even homotopy invariants, even though the left-hand side
of the isomorphism in (1) is obviously defined in terms of the
differentiable structure on $M$.

However, \deRham's theorem asserts still more, namely that the
isomorphism $R^p(\halfneg M)\Totto\halfneg H^p(\halfneg M)$ is
realized through the \textit{\textsf{integration}\/} of differential
forms. We assume that the cells $\sigma_p$ are differentiable and use
them as domains of integration; thus, for example,
$\int_{\sigma_1}\alpha^1$, $\alpha^1\in\Omega^1$, is an ordinary line
integral. For $\alpha^p\in\Omega^p$, $\int_{\sigma_p}\alpha^p$
extends, through linearity, to all chains $c_p\in C_p$;
$f^p(c_p)=\int_{c_p}\alpha^p$ defines a linear form $f^p\in C^p$, and
we obtain linear mappings
\[ \Phi:\Omega^p\Totto C^p
   \textrm{,}\quad\quad p\ =\ 0, 1, ..., n \textrm{.}
\]
For these, the equation
\begin{equation}
  \Phi(\dee{\alpha^p})\ =\ \delta\Phi(\alpha^p)
\end{equation}
holds. This is nothing other than Stokes's theorem:
\[ \singleIntegral{\dee{\alpha^p}}{}{c_{p+1}}{}\ =\ 
   \singleIntegral{\alpha^p}{}{\partial c_{p+1}}{}\ =\ f^p(\partial
   c_{p+1})\ =\ \delta f^p(c_{p+1})
\]
for all $c_{p+1}\in C_{p+1}$.

From (2), we see that $\Phi$ maps the groups $\mathring\Omega^p$ and
$\dee{\Omega^{p-1}}$ into $\mathring{C}^p$ and $\delta C^{p-1}$,
respectively, and hence gives rise to a linear mapping of the quotient
groups $R^p(M)$ into $H^p(M)$---and \deRham's theorem asserts that
this is the isomorphism in (1).

This is not exactly the original version given in the
\textit{\textsf{th\`ese}\/}. In order to arrive at the original form,
one has to consider the homology groups $H_p(M)$, that is, the
quotient groups $\mathring{C}_p/\partial C_{p+1}$, where
$\mathring{C}_p$ is the kernel of $\partial: C_p\Totto
C_{p-1}$. Chains $c_p$ such that $\partial c_p=0$ are called
``cycles,'' and those such that $c_p=\partial c_{p+1}$,
``boundaries''; they have a clear geometric significance. Cycles are
closed domains of integration, and an element of $H_p(M)$ is a class
of $p$-cycles differing from one another ``only'' by boundaries.

Now if $f^p\in\mathring{C}^p$, $\delta f^p=0$, then $f^p(\partial
c_{p+1})= \partial f^p(c_{p+1})=0$, that is, $f^p$ is $0$ on $\partial
C_{p+1}$. If we consider $f^p$ only on the cycles of $\mathring{C}_p$,
then we thereby obtain a linear form on
$H_p(M)=\mathring{C}_p/\partial C_{p+1}$; this is identically $0$ if
$f^p=\delta g^{p-1}$. Conversely, we obtain all real linear forms
$H_p(M)\Totto\mathbb{R}$ in this way (because a linear form on
$\mathring{C}_p$ can be extended to the entire vector space $C_p$ in
an arbitrary manner). Therefore
\[ H^p(M)\ =\ \Hom(H_p(M), \mathbb{R})
   \textrm{,}
\]
where $\Hom(V,\mathbb{R})$ denotes the vector space of all real linear
forms on $V$, that is, the vector space dual to $V$. The mapping $\Phi$ can
now be interpreted as $R^p(M)\Totto\Hom(H_p(M),\mathbb{R})$. It associates
a closed differential form $\alpha^p\in\Omega^p$ to a linear form that is
defined on the cycles and is equal to $0$ on the boundaries; its values are
called the \textit{\textsf{periods}\/} of $\alpha_p$. That $\Phi$ is an
isomorphism means two things:
\begin{enumerate}
\item[\textsf{I.}] If all periods of the closed differential form
$\alpha^p$ are equal to $0$, then $\alpha^p$ is a differential
$\dee{\beta^{p-1}}$. In other words, $\Phi$ is injective.
\item[\textsf{II.}] Given specified periods
$\varphi\in\Hom(H_p(M),\mathbb{R})$, there exists a closed differential
form $\alpha^p$ that has exactly these periods. In other words, $\Phi$ is
surjective. 
\end{enumerate}\vspace{-1\parskip}
This was the form in which \deRham\ worked out the proof: at the time,
this was a monumental achievement! Upon closer analysis it reveals
itself to be very similar to subsequent proofs, which were rendered
simple and transparent thanks to the techniques of sheaf theory (a
formalization of the passage from the local to the global).

In addition to Theorems \textsf{I} and \textsf{II}, the
\textit{\textsf{th\`ese}\/} also contains yet another, third
theorem. It deals with the (exterior) product $\alpha^p\wedge\beta^q$
of differential forms, for which the equation
\[ \dee{(\alpha^p\wedge\beta^q)}\ =\ \dee{\alpha^p}\wedge\beta^q + (-1)^p
   \alpha^p\wedge\dee{\beta^q}
\]
holds; it makes the direct sum $R^*(M)=\bigoplus R^p(M)$ into a graded
ring. Likewise, $H^*(M)=\bigoplus H^p(M)$ is a graded ring by means of
the \textit{\textsf{Alexander}\/} multiplication. We then have that
$\Phi$ is product-preserving, and thus a ring isomorphism $R^*(M)\cong
H^*(M)$. Of course, \deRham\ did not have the Alexander product at his
disposal; however, according to Poincar\'e duality, $H^p(M)\cong
H_{n-p}(M)$, and thus to each closed form $\alpha^p$ there corresponds
a cycle and hence a homology class $\Psi(\alpha^p)\in H_{n-p}(M)$. De
Rham proved that
\begin{enumerate}
\item[\textsf{III.}] \hfill $\Psi(\alpha^p\wedge\beta^q)\ =\
 S(\Psi(\alpha^p), \Psi(\beta^q))$, \hfill\ \break
\end{enumerate}\vspace{-2\parskip}
where $S$ denotes the ``intersection'' of the cycles $c_{n-p}=
\Psi(\alpha^p)$ and $d_{n-q}=\Psi(\beta^q)$. $S(c_{n-p}, d_{n-q})$ lies in
dimension $n-(p+q)$ and is Poincar\'e-dual to the Alexander product of $f^p
=\Phi(\alpha^p)$ and $g^q=\Phi(\beta^q)$. This theorem expresses a still
closer relationship between the algebra of differential forms and the
topology of the manifold.

\vspace{\parskip}
\textbf{\large Outlook}

De~Rham's theorems, which have only been outlined here, simultaneously
signified the end of a line of development and a new beginning.  On
one hand, the theory of differential forms growing out of analysis and
differential geometry and, on the other hand, the emergence of the 
concept of homology in combinatorial topology had led \'Elie
Cartan---and perhaps earlier even Poincar\'e---to conjecture that
these theorems must be true.  But the proof had to wait until \deRham\ 
succeeded in providing the necessary link from the local to the
global.  Not only the results, but also the methods would soon lead to
a wellspring of new ideas, which would later be further developed in
the works of Weil, Henri Cartan, Serre, Grothendieck, Sullivan, and
others; these would be concerned not only with manifolds (real, open,
complex, algebraic), but also with much more general cell complexes.

Georges \deRham's manner of presentation, both verbally and in his
writings, was crystal clear and ``concrete'' in the best sense of the
word. Any kind of ``abstract generalized nonsense,'' if not rooted in
some way in a concrete problem, was foreign to his thinking. Yet
according to my personal recollection, he also appreciated the value
and the expressive power of generalizations, especially of structural
concepts that were abstracted from mathematical substance. He was
thoroughly positive toward the efforts of the young Bourbaki group, to
which he was linked by genuine friendships. And he would without a
doubt have taken great joy in today's mathematics, with its numerous
surprising connections between seemingly disparate areas and its
synthesis of concrete problems and abstract general methods.

\vfill\pagebreak
{\small Professor Beno Eckmann,\\
Forschungsinstitut f\"ur Mathematik\\
ETH-Zentrum,\\
CH-8092 Z\"urich}

\end{document}